\documentclass[12pt]{amsart}
\usepackage{amssymb, mathrsfs,amsmath}
\usepackage{latexsym}
\usepackage{amsfonts}
\usepackage{shadow}
\newcommand{\longhookrightarrow}
{\ensuremath{\lhook\joinrel\relbar\joinrel\relbar\joinrel\rightarrow}}
\date{}
\author[D. Alpay]{Daniel Alpay}
\author[P. Jorgensen]{Palle Jorgensen}
\address{(DA) Department of Mathematics \newline
Ben Gurion University of the Negev \newline P.O.B. 653, \newline
Be'er Sheva 84105, \newline ISRAEL} \email{dany@math.bgu.ac.il}
\address{(PJ)
Department of Mathematics\newline 14 MLH \newline The University
of Iowa Iowa City,\newline IA 52242-1419 USA}
\email{jorgen@math.uiowa.edu}
\author[D. Levanony]{David Levanony}
\address{(DL) Department of Electrical Engineering \newline
Ben Gurion University of the Negev \newline P.O.B. 653, \newline
Be'er Sheva 84105, \newline ISRAEL }
\email{levanony@ee.bgu.ac.il}
\thanks{D. Alpay thanks the Earl Katz family for endowing the chair which
supported his research. The research of the authors was supported
in part by the Israel Science Foundation grant 1023/07. The work
was done while the second named author visited Department of
Mathematics, Ben Gurion University of the Negev, supported by a
BGU distinguished visiting scientist program. Support and
hospitality is much appreciated. We acknowledge discussions with
colleagues there, and in the US, Dorin Dutkay, Myung-Sin Song,
and Erin Pearse.}


\title[Stochastic processes with fractional measures]
{A class of Gaussian processes with fractional spectral measures}
\begin{document}
\parindent 0cm
\newtheorem{Pa}{Paper}[section]
\newtheorem{Tm}[Pa]{{\bf Theorem}}
\newtheorem{La}[Pa]{{\bf Lemma}}
\newtheorem{Cy}[Pa]{{\bf Corollary}}
\newtheorem{Rk}[Pa]{{\bf Remark}}
\newtheorem{Pn}[Pa]{{\bf Proposition}}
\newtheorem{Dn}[Pa]{{\bf Definition}}
\newtheorem{Ex}[Pa]{{\bf Example}}
\numberwithin{equation}{section}
\def\L{\mathbf L}
\def\R{\mathbb R}
\def\N{\mathbb N}
\def\C{\mathbb C}
\def\s{\mathscr S}
\def\ss{\mathscr S^\prime}
\def\sr{\mathscr S(\R)}
\def\ssr{\mathscr S'(\R)}
\def\(s){\mathscr S(\R^2)}
\def\F{\mathcal F}
\def\P{\mathcal P}
\def\W{\mathcal W}
\def\w{\omega}
\def\Dom{{\rm dom}~(T_m)}
\def\Doms{{\rm dom}~(T_m^*)}
\def\Def{\stackrel{{\rm def.}}{=}}
\begin{abstract}
We study a family of stationary increment
Gaussian processes, indexed by time.
These processes are determined by certain measures $\sigma$
(generalized spectral measures), and our focus here is on the case
when the measure $\sigma$ is a singular measure. We characterize
the processes arising from when $\sigma$ is in one of the classes
of affine self-similar measures. Our analysis makes use of
Kondratiev-white noise spaces.  With the use of {\it a priori}
estimates and the Wick calculus, we extend and sharpen (see
Theorem \ref{bastille}) earlier computations of Ito stochastic
integration developed for the special case of stationary
increment processes having absolutely continuous measures. We
further obtain an associated Ito formula (see Theorem
\ref{oberkampf}).
\end{abstract}
\maketitle
\tableofcontents
\section{Introduction}
\setcounter{equation}{0}
%
There are two ways of looking at stochastic processes, i.e.,
random variables indexed by a continuous parameter (for example
time): (i) One starts with a probability space, i.e., a sample
space, a set  $\Omega$ with a sigma algebra $\mathcal B$ of
subsets, and a probability measure $P$ on  $(\Omega, \mathcal B)$,
and a system of random variables $\left\{X(t)\right\}$ on the
$(\Omega, \mathcal B, P)$. From this one may then compute quantities such
as means, variances, co-variances, moments, etc, and then derive
important {\bf spectral data} . These in turn are used in
various applications, such as in solving stochastic
differential equations. Here we are concerned with the other
direction: (ii) Given some {\it  a priori} spectral data, how do
we construct a suitable probability space $(\Omega, \mathcal B,
P)$ and an associated process $\left\{X(t)\right\}$  such that the
prescribed spectral data is recovered from the constructed
process? In other words, this is a version of an {\it inverse
spectral problem}. For a number of reasons, it is useful in the
study of the inverse problem to focus on the case of Gaussian
processes.\\

A zero mean Gaussian process $\left\{X(t)\right\}$ on a
probability space is said to be {\it stationary increment} if the
mean-square expectation of the increment $X(t)-X(s)$ is a
function only of the time difference $t-s$. Then there is a
measure $\sigma$ such that the {\it covariance function} of such
a process is of the form
\begin{equation}
\label{Ksigma} E[X(t)X(s)^*]=K_\sigma(t,s)=\int_{\mathbb
R}\chi_t(u)\chi_s(u)^*d\sigma(u),\quad t,s\in\mathbb R,
\end{equation}
where $E$ is expectation, and we have set
\begin{equation}
\chi_t(u)=\frac{e^{itu}-1}{u}.
\end{equation}
The positive measure $d
\sigma$ is called the {\it spectral
measure}, and is subject to the restriction
\begin{equation}
\int_{\mathbb R}\frac{d\sigma(u)}{u^2+1}<\infty.
\label{eq:cond1}
\end{equation}
The covariance function $K_\sigma(t,s)$ can be rewritten as
\[
K_\sigma(t,s)=r(t)+r(s)^*-r(t-s),
\]
where
\begin{equation}
\label{kipour_2007}
r(t)=-\int_{\R}
\Big\{e^{itu}-1-\frac{itu}{u^2+1}
\Big\}\frac{d\sigma(u)}{u^2},
\end{equation}
When $\sigma$ is even, $r$ is real and takes the simpler form

\begin{equation}
\label{rreal}
r(t)=\int_{\mathbb R} \dfrac{1-\cos(tu)}{u^2}d\sigma(u).
\end{equation}

We note that some authors call spectral measure instead the
measure $u^2d\sigma(u)$ rather than the
measure $d\sigma(u)$. See \cite[p. 25 (7)]{Lifshits95}.\\

The literature contains a number of papers dealing with these
processes, but our treatment here goes beyond this, offering two
novelties: the inverse problem (see above), and an operator
theory of singular measures.  Both are motivated by the need to
deal with families of singular measures $\sigma$ (see
\eqref{Ksigma} through \eqref{kipour_2007}). Our focus is on
families of purely singular measures $\sigma$ with an intrinsic
{\sl spatial} selfsimilarity, typically with Cantor support, and
with fractional scaling (and Hausdorff) dimension; see Section
\ref{sec:self} below; these are measures with affine
selfsimilarity. Note that this notion is different from
self-similarity in the {\sl time}-variable; the latter case
includes fractional Brownian motion (fBm), studied in e.g.,
\cite{ptrf,MR2414165,aal2,MR1408433,JoSo09b}. For the latter
(fBm), it is known that the corresponding one-parameter family of
measures ${\sigma}$ consists of a scale of absolutely continuous
measures.\\

The derivative of a stationary increment process is a (possibly generalized)
stationary process, with covariance function
\[
\widehat{\sigma}(t-s),
\]
where $\widehat{\sigma}$ denotes the Fourier transform, possibly
in the sense of distributions, of $\sigma$. For a function $f$,
recall the Fourier transform
\[
\widehat{f}(u)=\int_{\mathbb R}e^{iux}f(x)dx.
\]

We note that second order stationary processes can be studied with
the use of the theory of Hilbert spaces and of unitary
one-parameter groups of operators in Hilbert space. One may then
invoke the Stone-von Neumann spectral theorem, the spectral
representation theorem, and a detailed multiplicity theory to
study these processes. See for instance \cite{KR97}.\\

An important role in the theory is played by the space $\mathcal
M(\sigma)$ of functions in $\mathbf L_2(\mathbb R,dx)$ such that
\[
\int_{\mathbb R}|\widehat{f}(u)|^2d\sigma(u)<\infty.
\]
%
This space contains in particular the Schwartz space.\\

In the paper \cite{aal2}, see also \cite{MR2414165}, the case
where $\sigma(u)$ is absolutely continuous with respect to
Lebesgue measure, i.e., $d\sigma(u)=m(u)du$ (where the
Radon-Nikodym derivative $m$ satisfies moreover some growth
conditions) was considered. The study of \cite{aal2} included in
particular the case of the Brownian motion and of the fractional
Brownian motion. A key role in that paper was played by the (in
general unbounded) operator $T_m$ on ${\mathbf L}_2(\mathbb R,dx)$
defined by
\begin{equation}
\label{tm}
\widehat{T_mf}=\sqrt{m}\widehat{f}.
\end{equation}
So $T_m$ is a convolution operator in $\mathbf L_2(\mathbb
R,dx)$, i.e.,
\[
T_m f=(\sqrt{m})^\vee\star f,
\]
with $\vee$ denoting the inverse Fourier transform, in the sense
of distributions.\\

In this paper we focus on the case when the spectral measure is
an {\it affine iterated function-system measure} (AIFSs). Among
the AIFS-measures there is a subfamily which admits an
orthonormal family of Fourier frequencies. These are lacunary
Fourier series studied first in a paper by one of the authors and
Steen Pedersen in 1998, see \cite{JoPe98}. A lacunary Fourier
series is one in which there are large gaps between consecutive
nonzero coefficients. AIFS-measures may be visualized as fractals
in the small, while their Fourier expansions as dual fractals in
the large. The spectral measure of such a process
$\left\{X(t)\right\}$ is important as it enters in a rigorous
formulation of an associated Ito formula for functions $f(X(t))$ of the
process.\\

The main results of the paper may be summarized as follows: We
construct a densely defined operator $Q_\sigma$ from $\mathbf
L_2(d\sigma)$ into $\mathbf L_2(\mathbb R,dx)$ such that

\begin{equation}
\label{magic2}
\int_{\mathbb R}\chi_t(u)\chi_s(u)^*d\sigma(u)=
\langle
Q_\sigma({1_{[0,t]}}),Q_\sigma({1_{[0,s]}})\rangle_{\mathbf
L_2(\mathbb R,dx)}
\end{equation}
This operator $Q_\sigma$ is the counterpart of the operator $T_m$
defined in \eqref{tm} and introduced in \cite{aal2}. We denote by
$f\mapsto\widetilde{f}$ the natural isometric imbedding of
${\mathbf L}_2(\mathbb R,dx)$ into the white noise space; See
Section \ref{2} for details. The stochastic process
$\left\{X_\sigma(t)\right\}$ defined by
\[
X_\sigma(t)=\widetilde{Q_\sigma({1_{[0,t]}})},\quad t\in\mathbb R,
\]
has covariance function
\[
E[X_\sigma(t)X_\sigma(s)^*]= \int_{\mathbb
R}\chi_t(u)\chi_s(u)^*d\sigma(u)
\]
Following \cite{Lifshits95}, the  measure $\sigma$ in this
expression will be called the {\it spectral measure} of the
process. Its intuitive meaning is that of  "spectral densities",
not to be confused with "power spectral measure" traditionally
used for the much more restrictive family of stochastic process,
the stationary processes. In the case of stationary processes,
and when the power spectral measure is absolutely continuous with
respect to Lebesgue measure, one speaks of power spectral density
(psd). It is then the Fourier transform of the covariance
function, a function of a single variable, namely, the time
difference. When the process can be differentiated, its
derivative is stationary and $\sigma$ is absolutely continuous
with respect to Lebesgue measure, its derivative
is the psd of the derivative process.\\

We show that $X_\sigma(t)$ admits a derivative
$(X_\sigma(t)^\prime)\stackrel{\rm def.}{=}W_\sigma(t)$ in the
white noise space, which is moreover continuous in the white noise
space norm. Furthermore
\[
E[W_\sigma(t)(W_\sigma(s))^*]=\widehat{\sigma}(t-s).
\]
It is found that both the processes $\left\{W_\sigma(t)\right\}$ and
$\left\{X_\sigma(t)\right\}$
are in the white noise space. We define a stochastic integral with
respect to $X_\sigma$, and prove results similar to those of \cite{ptrf},
but for a different class of processes studied here.\\

The outline of the paper is as follows. The paper consists of nine
sections besides the introduction. The first three small sections
are of a review nature. In Section \ref{sec:self} we present some
material on measures with affine selfsimilarity. In Section
\ref{3} we recall some properties of the associated $\mathbf
L_2(d\sigma)$ spaces. Hida's white noise space theory is based on a
Hilbert space, and plays an important role in our work. Its main
features are listed in Section \ref{2}. In Sections 4-9 we
develop the new results of the paper. In Section \ref{secQ} we
construct an operator for which \eqref{magic2} holds. The
corresponding process $X_\sigma$ and its derivative are
constructed in Section \ref{secXW}. The associated stochastic
integral and a Ito formula are considered in Sections
\ref{wickito} and \ref{itoform} respectively. To contrast with
the measures considered here, two examples of stationary
increments Gaussian processes with measures with unbounded
support are presented in Section \ref{two_ex}. In Section
\ref{QG} we briefly consider the case of a general measure
$\sigma$. The last section is devoted to various concluding
remarks.

\section{Measures with affine selfsimilarity}
\setcounter{equation}{0}
\label{sec:self}
In understanding
processes $\left\{X(t)\right\}$ with stationary increments, one
must look at the variety possibilities for measures $\sigma$,
representing spectral measures, in the sense outlined above. Each
measure-type for $\sigma$ entails properties of the associated
Ito formulas for $\left\{X(t)\right\}$ as it enters
into stochastic integration formulas.\\

While earlier literature on stationary-increment processes
has been focused on the case when $\sigma$ was assumed to be
absolutely continuous with respect to Lebesgue measure, or
perhaps the case when it is singular but purely atomic; in this
section we will focus instead on a quite different family of
measures: purely singular and non-atomic.
They share the following four features:\\

$(i)$ they are given by explicit recursive formulas;\\
$(ii)$ they posses an intrinsic affine selfsimilarity; see
\eqref{eq:AIFS},\\
$(iii)$ they admit a harmonic analysis based
on a lacunary
Fourier expansion; see \eqref{eq:lambdam}, and finally,\\
$(iv)$ the Fourier transform of $\sigma$ admits an explicit
infinite-product formula; see \eqref{sigmarho}.

\begin{Dn}
A Borel probability measure $\sigma$ on $\mathbb R$ is said to be an
affine iterated function system measure (AIFS) if there is a
finite family $\mathcal F$ of (usually contractive) affine
transformations on $\mathbb R$ such that
\begin{equation}
\label{eq:AIFS}
\sigma=\frac{1}{{\rm card}~{\mathcal
F}}\sum_{\tau\in\mathcal F} \sigma\circ\tau^{-1}
\end{equation}
holds i.e.
\[
\int f(x)d\sigma(x)=\frac{1}{{\rm card}~{\mathcal F}}\sum_{\tau\in
\mathcal F}\int f(\tau(x))d\sigma(x),
\]
for all bounded continuous functions on $\mathbb R$.
\end{Dn}
The simplest examples are Bernouilli convolutions. Then ${\rm
card}~{\mathcal F}=2$ and there exists some fixed $\rho>0$ such that
\begin{equation}
\label{rhopm}
\tau_{\pm}(x)=\rho(x\pm1).
\end{equation}
In that case, the measure $d\sigma_\rho$ satisfying
\eqref{eq:AIFS}, has Fourier transform of the form
\begin{equation}
\label{sigmarho}
\widehat{\sigma_\rho}(t)=\prod_{k=1}^\infty
\cos(\rho^kt).
\end{equation}
Note that the function
\[
\prod_{k=1}^\infty \cos(\rho^k(t-s))
\]
is positive definite on the real line, since each of the
functions in that product
\[
\cos(\rho^k(t-s))=
\cos(\rho^kt)\cos(\rho^ks)+\sin(\rho^kt)\sin(\rho^ks)
\]
is positive definite, and one can obtain $\sigma$ from
Bochner's theorem.\\

Cases with $\rho$ of the form
\[
\rho=\frac{1}{2m},\quad m=2,3,4,\ldots
\]
will be of special interest here. Fix $m$ and let $\sigma_m$
(that is with $\rho=\frac{1}{2m}$) be the corresponding
Bernouilli measure. For $t\in\mathbb R$, set
\[
e_t(u)=e^{ i tu}.
\]
Let
\begin{equation}
\label{eq:lambdam}
\Lambda_m=2\pi\left\{\sum_{0}^N
b_j(2m)^j,\,\,{\rm where}\,\,N\in\mathbb N\,\, {\rm and}\,\ b_j\in
\left\{0,\frac{m}{2}\right\}\right\}.
\end{equation}
For instance,
\[
\begin{split}
\Lambda_2&=2\pi\left\{0,1,4,5,16,17,20,\ldots\right\}\\
\Lambda_3&=2\pi\left\{0,\frac{3}{2},9,\frac{21}{2},18,\ldots\right\}\,\,\,
{\rm and}\\
\Lambda_4&=2\pi\left\{0,2,16,18,128,130,\ldots\right\}.
\end{split}
\]

It is known that the set $\Lambda_m$ makes
\[
\left\{\, e_\lambda\, |\,\, \lambda\in\Lambda_m\right\}
\]
into an orthonormal basis (ONB) in $\mathbf L_2(d\sigma_m)$; we
say that $(\sigma_m,\Lambda_m)$ is a {\it spectral pair}. Let us
formalize this notion:

\begin{Dn}
A Borel finite measure $\sigma$ on $\mathbb R$ is said to have a
spectrum $\Lambda\subset\mathbb R$ if $\Lambda$ is a discrete set
and the set $\left\{e_\lambda,\,\, \lambda\in\Lambda\right\}$ is
an orthonormal basis in $\mathbf L_2(d\sigma)$. Then
$(\sigma,\Lambda)$ is called a spectral pair.
\end{Dn}

By the Fourier basis property mentioned above, we refer to
the presence of a Fourier orthonormal basis (ONB) in the Hilbert
space $\mathbf L_2(d\sigma)$; and our discussion below is restricted
to the case when $\sigma$ is assumed to be a finite
measure. The study of these singular measures was initiated by
one of the authors in collaboration with co-authors, see \cite{
JoPe98, Str98, JoPe95, JoPe94, JoPe93a, JoPe93b, JoPe92, DuJo06,
JKS07, DDJ09, DuJo09, DuJo08}.
The Fourier expansion in $\mathbf L_2(d\sigma)$
for fractal measures $\sigma$
  differs from standard Fourier
series (for periodic functions) in that the fractal Fourier
expansion is local, much like wavelet expansions; see \cite{Str00}
for details. While the family of these singular measures is
extensive, we found it helpful to focus our discussion below
on one of the simplest cases, the first occurring in the
literature, see \cite{JoPe98}. It has Hausdorff dimension= scaling
dimension= 1/2, and its support is a Cantor-subset of the real
axis.\\

It is proved in \cite{JoPe98} that the support of $d\sigma_m$ is
inside the closed interval $[-1/2,1/2]$, and has Lebesgue measure
$0$. We note however, there are also spectral pairs
$(\sigma_m,\Delta_m)$ where the measure $d\sigma_m$  is not
compactly supported.

\section{The spaces $\mathbf L_2(d\sigma)$}
\setcounter{equation}{0} \label{3}
For later use, we review two results from
\cite{DuJo06,DDJ09,DuJo08,DuJo09} and \cite{JoPe98}. Recall that,
for $t\in\mathbb R$, $e_t(u)=e^{ i tu}$. In general one does not
assume that $\sigma$ has compact support. When the support of
$\sigma$ is compact, it has a well defined Fourier transform,
which is an entire function and not merely a distribution. We
have the following result, proved in \cite{JoPe94,JoPe95,DuJo09}.

\begin{Tm}
Let $\sigma$ be a finite positive Borel measure and let
$\Lambda\subset\mathbb R$ be a discrete set. Then, $(\sigma,
\Lambda)$ is a spectral pair if and only if
\begin{equation}
\label{eq:equiv1}
\sum_{\lambda\in\Lambda}|\widehat{\sigma}(t-\lambda)|^2=1,\quad\forall
t\in \mathbb R.
\end{equation}
\end{Tm}

\begin{La}
Let $t,s\in\mathbb R$. It holds that
\begin{equation}
\|e_t-e_s\|_{{\mathbf L_2(d\sigma)}}\le K|t-s| \label{ineq:ets}
\end{equation}
where
\[
K=\int_{[-\frac{1}{2},\frac{1}{2}]}u^2d\sigma(u).
\]
\end{La}

{\bf Proof:} We have
\[
\begin{split}
\|e_t-e_s\|_{{\mathbf
L_2(d\sigma)}}^2&=\int_{[-1/2,1/2]}|e^{itu}-e^{isu}|^2d\sigma(u)\\
&=\int_{[-1/2,1/2]}|1-e^{iu(t-s)}|^2d\sigma(u)\\
&=\left(\int_{[-1/2,1/2]}4u^2\frac{\sin^2\left(\frac{u(t-s)}{2}\right)
}{\left(\frac{u(t-s)}{2}\right)^2}d\sigma(u)\right)\cdot\frac{(t-s)^2}{4}\\
&\le\left(\int_{[-1/2,1/2]}u^2d\sigma(u)\right)\cdot(t-s)^2.
\end{split}
\]
\mbox{}\qed\mbox{}\\

As already mentioned, the measures $\sigma$ we consider are such
that an orthonormal basis of $\mathbf L_2(d\sigma)$ is of the form
\[
e^{i \lambda_n u},n=0,1,\ldots,
\]
where $\lambda_n\in \pi\mathbb  N_0$ for all $n\in\mathbb N_0$.

\section{A brief survey of white noise space analysis}
\label{2}
\setcounter{equation}{0}
In this section we present some technical details required in the
subsequent sections, taken from Hida's white noise space theory. We
refer the reader to \cite{MR1244577}, \cite{MR2444857},
\cite{MR1408433} for more information. The facts reviewed here
are essential for our analysis of certain stochastic integrals
(Section \ref{wickito}), and our Ito formula (Section
\ref{itoform}). While convergence questions for stochastic
integrals traditionally involve integration in probability spaces
of paths, in our approach, the sample space will instead be a
space of tempered distributions $\mathcal S^\prime$ derived from
a Gelfand triple construction; but there is a second powerful
tool involved, a completion called the Kondratiev-Wick algebra.
We briefly explain the justification for this approach below.\\

The second system of duality spaces are called Kondratiev spaces,
see Section \ref{2} for details. Further, there is a particular
Kondratiev space, endowed with a product, the Wick product and
an algebra under this product. It serves as a powerful tool in
building stochastic integrals because, as we show, the stochastic
integral takes place in the Kondratiev-Wick algebra; and we can
establish convergence there; see Theorem \ref{bastille}. Moreover
(see Theorem \ref{oberkampf}), the stochastic integration making
up
our Ito formula lives again in the Kondratiev-Wick algebra.\\

Let $\mathcal S$ denote the Schwartz space of {\it real-valued}
$C^\infty(\mathbb R)$ functions such that
\[
\forall p,q\in\mathbb
N_0,\quad\lim_{x\rightarrow\pm\infty}x^pf^{(q)}(x)=0
\]
For $s\in{\mathcal S}$, let $\|s\|$ denote its ${\mathbf
L}
_2({\mathbb R})$ norm. The function
\[
K(s_1-s_2)=e^{-\frac{\|s_1-s_2\|^2}{2}}
\]
is positive definite for $s_1,s_2$
running in ${\mathcal S}$. By the Bochner-Minlos theorem (see \cite{MR0154317},
\cite[Th\'eor\`eme 3, p. 311]{MR35:7123}), there exists a
probability measure $P$ on ${\mathcal S}^\prime$ such that
\begin{equation}
\label{KS}
K(s)=\int_{{\mathcal S}^\prime}e^{-i\langle
s^\prime,s\rangle}dP(s^\prime),
\end{equation}
where $\langle s^\prime,s\rangle$ denotes the duality
between ${\mathcal S}$ and ${\mathcal S}^\prime$.
Henceforth, we set
$\Omega=\mathcal S^\prime$.
The real Hilbert space $\mathcal
W={\mathbf L}_2(\Omega, {\mathscr F},dP)$, where ${\mathscr F}$
is the Borelian $\sigma$-algebra, is called the
white noise space.\\

For $s\in\mathcal S$ and $\omega\in\Omega$ we set
\begin{equation}
\widetilde{s}(\w)=\langle \w,s\rangle
\label{ftilde}
\end{equation}
From \eqref{KS} follows that the map $s\mapsto \widetilde{s}$ is an isometry
from $\mathcal S$ endowed with the $\mathbf L_2(\mathbb R,dx)$ norm into
$\mathcal W$. This isometry extends to all of $\mathbf L_2(\mathbb R,dx)$, and
we will denote the extension by the same symbol.\\

We now present an orthogonal basis of $\mathcal W$. We set $\ell$
to be the space of sequences $(\alpha_1,\alpha_2,\ldots)$, whose
entries are in
\[
\mathbb N_0=\left\{ 0,1,2,3,\ldots\right\},
\]
where $\alpha_k\not =0$ for only a finite number of indices $k$.
Furthermore, we denote by $\mathbf h_0,\mathbf h_1,\ldots$ the
Hermite polynomials. The functions
\[
H_\alpha=H_\alpha(\omega)
=\prod_{k=1}^\infty {\mathbf
h}_{\alpha_k}(\widetilde{h_k}(\w)),\quad \alpha\in\ell,
\]
form an orthogonal base of the white noise space
(the $\omega$-dependence will be omitted throughout, unless specifically
required). Furthermore, one has
\begin{equation}
\label{eq:fock1}
\|H_\alpha\|_{\mathcal W}^2=\alpha!,
\end{equation}
where we have used the multi-index notation
\[
\alpha!=\alpha_1!\alpha_2!\cdots,
\]
The Wick product $\lozenge$ in ${\mathcal W}$ is defined by the
formula
\[
H_\alpha\lozenge H_\beta=H_{\alpha+\beta},\quad \alpha,
\beta\in\ell,
\]
on the basis $(H_\alpha)_{\alpha\in\ell}$, and is extended by
linearity to $\mathcal W$ as
\begin{equation}
\label{wertyu} F\lozenge
G=\sum_{\gamma\in\ell}(\sum_{\alpha+\beta=\gamma}f_\alpha
g_\beta)H_\gamma,
\end{equation}
where $F=\sum_{\alpha\in\ell}f_\alpha H_\alpha$ and
$G=\sum_{\alpha\in\ell}g_\alpha H_\alpha$. See \cite[Definition
2.4.1, p. 39]{MR1408433}. The Wick product $F\lozenge G$ reduces
to multiplication by a constant when one of the elements $F$ or
$G$ is non random. The Wick product is not everywhere defined in
$\mathcal W$, and one may remedy this by viewing $\mathcal W$ as
the middle part of a Gelfand triple. The first element in the
triple is the Kondratiev space $S_1$ of stochastic test
functions, defined as the intersection of the Hilbert spaces
$\mathcal H_k$, $k=1,2,\ldots$, of series
$f=\sum_{\alpha\in\ell}f_\alpha H_\alpha$ such that
\begin{equation}
\label{michelle111}
|\|f\||^2_{k}\stackrel{\rm def.}{=}
\sum_{\alpha\in\ell}(\alpha!)^2|f_\alpha|^2 (2{\mathbb
N})^{k\alpha}<\infty.
\end{equation}

The third element in the Gelfand triple is the Kondratiev space
$S_{-1}$ of stochastic distributions. It is
a nuclear space, and is defined as the inductive limit of the
increasing family of Hilbert spaces ${\mathcal H}_{k}^\prime,
k=1,2,\ldots$ of formal series $\sum_{\alpha\in\ell}f_\alpha
H_\alpha$ such that
\begin{equation}
\label{michelle}
\|f\|^2_{k}\stackrel{\rm def.}{=} \sum_{\alpha\in\ell}|f_\alpha|^2
(2{\mathbb N})^{-k\alpha}<\infty,
\end{equation}where, for $\beta\in\ell$,
\[
(2{\mathbb N})^{\pm\beta}=2^{\pm\beta_1}(2\times 2)^{\pm\beta_2}
(2\times 3)^{\pm\beta_3}\cdots.
\]
See  \cite[\S 2.3, p. 28]{MR1408433}.\\

The Wick product is stable both in $S_1$ and in $S_{-1}$.
Moreover, V\r{a}ge's inequality (see \cite[Proposition 3.3.2,
  p. 118]{MR1408433}) makes precise the fact
  that $F\lozenge G\in
  S_{-1}$ for every choice of $F$ and $G$ in $S_{-1}$:
Let $l$ and $k$ be natural numbers such that $k>l+1$. Let $h\in
{\mathcal H}_{l}^\prime$ and $u\in {\mathcal H}_{k}^\prime$. Then,
\begin{equation}
\label{vage}
\|h\lozenge u\|_{k}\le A(k-l)\|h\|_{l}\|u\|_{k},
\end{equation}
where
\begin{equation}
\label{vage111} A(k-l)=\left(\sum_{\alpha\in\ell}(2{\mathbb
N})^{(l-k)\alpha}\right)^{1/2}.
\end{equation}

\section{The operator $Q_\sigma$}
\setcounter{equation}{0}
\label{secQ}
As we saw, the construction of a process $\left\{X(t)\right\}$
from a fixed spectral measure $\sigma$  depends on properties of a
certain operator in $\mathbf L_2(\mathbb R,dx)$. If $\sigma$  is
assumed absolutely continuous, with Radon-Nikodym derivative $m$,
then this operator $T_m $ was studied earlier and it is a
convolution operator with the square root of $m$. See \eqref{tm}
and \cite{aal2}. In this section we introduce the counterpart of
the operator $T_m$  in the present setting. Recall that
$h_0,h_1,\ldots$ denote the Hermite functions.  We define a
unitary map $W$ from $\mathbf L_2(d\sigma)$ onto $\mathbf
L_2(\mathbf R,dx)$ via the formula
\[
W(e_{\lambda_n})=h_n,\quad n=0,1,\ldots
\]
Let $M_u$ denote the operator of multiplication by the variable $u$
in $\mathbf L_2(d\sigma)$. The formula
\begin{equation}
\label{psiT}
T=WM_uW^*
\end{equation}
defines a bounded self-adjoint operator from $\mathbf L_2(\mathbf
R,dx)$ into itself. Furthermore, for $\psi\in\mathcal S$ we set
\[
Q_\sigma(\psi)=\widehat{\psi}(T)h_0.
\]

\begin{La}
Let $\psi\in\mathcal S$. Then, it holds that:
\begin{equation}
\label{Qpsi}
(Q_\sigma\psi)(x)=\sum_{n=0}^\infty\left(\int_{\mathbb
R}\widehat{\sigma}(y-\lambda_n)\psi(y)dy\right)h_n(x)
\end{equation}
\end{La}

{\bf Proof:} Using the functional calculus we have from \eqref{psiT}
that,
\[
\widehat{\psi}(T)=W\widehat{\psi}(M_u)W^*
\]
and hence
\[
\begin{split}
\widehat{\psi}(T)h_0&=W\widehat{\psi}(M_u)W^*h_0\\
&=W\widehat{\psi}(M_u)1\\
&=W\widehat{\psi}(u).
\end{split}
\]
Moreover, since $\widehat{\psi}\in{\mathbf L}_2(d\sigma)$ we have
\[
\begin{split}
\widehat{\psi}(u)&=\sum_{n=0}^\infty \langle
\widehat{\psi},e_{\lambda_n}
\rangle_{\mathbf L_2(d\sigma)}e_{\lambda_n}(u)\\
&=\sum_{n=0}^\infty \left(\int_{\mathbb R
}\widehat{\psi}(u)e^{-iu\lambda_n}
d\sigma(u)\right)e_{\lambda_n}(u),
\end{split}
\]
so that $W^*\widehat{\psi}$ is given by \eqref{Qpsi}.
\mbox{}\qed\mbox{}\\

The operator $Q_\sigma$ is typically an unbounded operator in the
Hilbert space $\mathbf L_2(\mathbb R,dx)$, but it is well defined
on a dense domain which consists of the Schwartz space $\mathcal
S\subset \mathbf L_2(\mathbb R,dx)$. These facts are elaborated
upon in the next theorem. In the statement, note that the
Fr\'echet topology of $\mathcal S$ is stronger than that of the
$\mathbf L_2(\mathbb R,dx)$-norm.

\begin{Tm}
Let $\psi\in\mathcal S$. Then, it holds that,
\begin{equation}
\label{normQ} \|Q_\sigma(\psi)\|^2_{\mathbf L_2(\mathbb
R,dx)}=\int_{\mathbb R}|\widehat{\psi}(u)|^2d\sigma(u).
\end{equation}
In particular, $Q_\sigma$ is a continuous operator from $\mathcal
S$ into $\mathbf L_2(\mathbb R,dx)$. More precisely,
\begin{equation}
\label{eq:new}
\|Q_\sigma\psi\|_{\mathbf L_2(\mathbb R,dx)} \le \sqrt{K}
\left(\left(\int_{\mathbb R}|\psi(x)|dx\right)^2+
\left(\int_{\mathbb R}|\psi^\prime(x)|dx\right)^2\right)^{1/2},
\end{equation}
\label{tm5.2} where
\begin{equation}
\label{K}
K=\int_{\mathbb R}\frac{d\sigma(u)}{1+u^2}<\infty.
\end{equation}
\end{Tm}

\begin{Rk}{\rm
For some stochastic processes it is important to relax condition
\eqref{K} to
\begin{equation}
\label{Kp}
K_p=\int_{\mathbb R}\frac{d\sigma(u)}{1+|u|^{2p}}<\infty
\end{equation}
for some $p\in\mathbb N$. In this case the estimate in
\eqref{eq:new} becomes
\begin{equation}
\label{eq:new1}
\|Q_\sigma\psi\|_{\mathbf L_2(\mathbb R,dx)} \le \sqrt{K_p}
\left(\left(\int_{\mathbb R}|\psi(x)|dx\right)^2+
\left(\int_{\mathbb R}|\psi^{(p)}(x)|dx\right)^2\right)^{1/2}.
\end{equation}
}
\end{Rk}

{\bf Proof of Theorem \ref{tm5.2}:} Let $\psi\in\mathcal S$. We
have:
\[
\begin{split}
\|Q_\sigma\psi\|_{\mathbf L_2(\mathbb R,dx)}^2&=\sum_{n=0}^\infty\big|
\int_{\mathbb R}\widehat{\sigma}(y-\lambda_n)\psi(y)dy \big|^2\\
&=\sum_{n=0}^\infty\big|\iint e^{iu(y-\lambda_n)}\psi(y)d\sigma(u)dy\big|^2\\
&=\sum_{n=0}^\infty\big|\int(\int \psi(y)e^{iuy}dy)e^{-iu\lambda_n}d\sigma(u)\big|^2\\
&=\sum_{n=0}^\infty\big|\int\widehat{\psi}(u)e^{-iu\lambda_n}d\sigma(u)\big|^2\\
&=\int_{\mathbb R}|\widehat{\psi}(u)|^2d\sigma(u),
\end{split}
\]
where we have used Fubini's theorem for the third equality, and
Parseval's equality for the last equality.\\

We now prove \eqref{eq:new}. This will prove the continuity of
$Q_\sigma$ from $\mathcal S$ endowed with its Fr\'echet topology
into $\mathbf L_2(\mathbb R,dx)$. For $\psi\in\mathcal S$ we have
\[
\begin{split}
\|Q_\sigma\psi\|^2_{\mathbf L_2(\mathbb R,dx)}&=\int_{\mathbb
R}|\widehat{\psi}(u)|^2d\sigma(u)\\
&=\int_{\mathbb
R}(1+u^2)|\widehat{\psi}(u)|^2\frac{d\sigma(u)}{1+u^2}\\
&\le K\max_{u\in\mathbb R}(1+u^2)|\widehat{\psi}(u)|^2 \\
&\le K\left(\int_{\mathbb R}|\psi\star\psi^\sharp|(x) dx+
\int_{\mathbb
R}|\psi^\prime\star(\psi^\sharp)^\prime|(x) dx\right)\\
&\le K\left(\left(\int_{\mathbb R}|\psi(x)|dx\right)^2+
\left(\int_{\mathbb R}|\psi^\prime(x)|dx\right)^2\right),
\end{split}
\]
where $\psi^\sharp(x)=\overline{\psi(-x)}$ and
$\psi^\prime=\frac{d\psi}{dx}$. Furthermore, we have used
\[
\widehat{\psi\star\psi^\sharp}(u)=|\widehat{\psi}(u)|^2.
\]
and
\[
\begin{split}
\int_{\mathbb R}|\psi\star\psi^\sharp|(x)
dx&\le\left(\int_{\mathbb R}|\psi|(x) dx\right)\left(\int_{\mathbb
R}|\psi^\sharp|(x) dx\right)\\
&=\left(\int_{\mathbb R}|\psi|(x) dx\right)^2.
\end{split}
\]

Since the expression on the right hand side in this estimate is
one of the Fr\'echet semi-norms of $\mathcal S$, the continuity
assertion in Theorem \ref{tm5.2} follows. The estimate
\eqref{eq:new} further gives an exact rate of continuity.
\mbox{}\qed\mbox{}\\

From equation \eqref{normQ} we can extend the domain of
definition of $Q_\sigma$ to a wider set, which in particular
include the functions $1_{[0,t]}$. This is explicited in the
following proposition. We remark that such a result may be extended to
more general measures $\sigma$'s, see Section \ref{QG}.

\begin{Pn}
\label{10t}
Let $f\in\mathbf L_2(d\sigma)$ be such that, for some
sequence $(s_n)_{n\in\mathbb N}$ of Schwartz functions,
\begin{equation}
\label{limsn}
\lim_{n\rightarrow\infty}|f-\widehat{s_n}|_\infty=0.
\end{equation}
Then the sequence $(Q_\sigma s_n)_{n\in\mathbb N}$ is a Cauchy sequence in
$\mathbf L_2(\mathbb R,dx)$. Its limit is the same for all
sequences which satisfy \eqref{limsn}, and will be denoted by
\[
Q_\sigma f\stackrel{\rm def.}{=}\lim_{n\rightarrow\infty}Q_\sigma s_n.
\]
\end{Pn}

{\bf Proof:} From \eqref{normQ} follows that for every
$\epsilon>0$ there is an $N\in\mathbb N$ such that
\[
n,m\ge N\Longrightarrow |\widehat{s_n}-\widehat{s_m}|_\infty\le
\epsilon.
\]
Thus for such $n$ and $m$
\begin{equation}
\label{eq_new}
\begin{split}
\|Q_\sigma s_n-Q_\sigma s_m\|^2_{\mathbf L_2(\mathbb
R,dx)}&=\int_{\mathbb R}
|\widehat{s_n}(u)-\widehat{s_m}(u)|^2d\sigma(u)\\
&\le \sigma(\mathbb
R)\cdot|\widehat{s_n}-\widehat{s_m}|^2_\infty\\
&\le \sigma(\mathbb R)\cdot\epsilon^2.
\end{split}
\end{equation}
Therefore, $\lim_{n\rightarrow\infty} Q_\sigma s_n$ in the norm of
$\mathbf L_2(\mathbb R,dx)$. Call this limit $q_1$, and assume
that, for another sequence $(t_n)_{n\in\mathbb N}$ satisfying
\eqref{limsn}, we obtain another limit, say $q_2$. Note that
\begin{equation}
\label{sumlim}
\lim_{n\rightarrow\infty}|\widehat{s_n}-\widehat{t_n}|_\infty\le
\lim_{n\rightarrow\infty}|\widehat{s_n}-f|_\infty+
\lim_{n\rightarrow\infty}|f-\widehat{t_n}|_\infty=0.
\end{equation}
Then,
\[
\begin{split}
\|q_1-q_2\|_{\mathbf L_2(\mathbb R,dx)}&\le\| q_1-Q_\sigma
s_n\|_{\mathbf L_2(\mathbb
R,dx)}+\\
&\hspace{5mm}+\|Q_\sigma s_n-Q_\sigma t_n\|_{\mathbf L_2(\mathbb
R,dx)}+
\|Q_\sigma t_n-q_2\|_{\mathbf L_2(\mathbb R,dx)}\\
&\le \| q_1-Q_\sigma s_n\|_{\mathbf L_2(\mathbb
R,dx)}+\\
&\hspace{5mm}+\sqrt{\sigma(\mathbb
R)}\cdot|\widehat{s_n}-\widehat{t_n}|_\infty+ \|Q_\sigma
t_n-q_2\|_{\mathbf L_2(\mathbb R,dx)},
\end{split}
\]
which goes to $0$ as $n\rightarrow\infty$ by definition of $q_1$
and $q_2$ and due to \eqref{sumlim}.\\

We now verify that the function
\[
\chi_t(u)=\frac{e^{itu}-1}{u}
\]
can be approximated in the supremum norm by Schwartz functions.
The function $\chi_t$ vanishes at infinity, and hence can be
approximated in the supremum norm by continuous functions with
compact support; see \cite[Theorem 3.17, p. 70]{Rudin}. These in
turn can be approximated by functions in $\mathcal S$, using
approximate identities, as for instance Step 6 in the proof of
Theorem 6.1 in \cite{ptrf}. We sketch the argument for
completeness. Let
\begin{equation}
\label{eq:sketch}
k_\epsilon(x)=\frac{1}{\sqrt{2\pi}\epsilon}\exp{\left(-\frac{x^2}{2\epsilon^2}
\right)}.
\end{equation}
$k_\epsilon$ is an $\mathcal N(0,\epsilon^2)$ density, and
therefore
\begin{equation}
\int_{\mathbb R} k_\epsilon(x)dx=1,
\label{int1}
\end{equation}
and, for every $r>0$
\begin{equation}
\label{int2}
\lim_{\epsilon\rightarrow
0}\int_{|x|>r}k_\epsilon(x)dx=0.
\end{equation}
Indeed, for $|x|>r>0$,
\[
\begin{split}
\frac{1}{\epsilon\sqrt{2\pi}}\int_{r}^\infty
e^{-\frac{x^2}{2\epsilon^2}}dx&=
\frac{1}{\epsilon\sqrt{2\pi}}\int_{r}^\infty \frac{x}{\epsilon^2}
e^{-\frac{x^2}{2\epsilon^2}}\frac{\epsilon^2}{x}dx\\
&\le \frac{\epsilon}{r\sqrt{2\pi}} \int_{r}^\infty
\frac{x}{\epsilon^2} e^{-\frac{x^2}{2\epsilon^2}}dx\\
&=\frac{\epsilon}{r\sqrt{2\pi}}e^{-\frac{r^2}{2\epsilon^2}}\\
&\longrightarrow 0\quad\mbox{\rm as}\quad \epsilon\rightarrow 0.
\end{split}
\]
Theses properties express the fact that $k_\epsilon$ is an
approximate identity. Applying \cite[Theorem 1.2.19, p.
25]{MR2463316} we see that, for every continuous function with
compact support,
\[
\lim_{\epsilon\rightarrow 0}\|k_\epsilon \ast f-f\|_\infty=0.
\]
To conclude, one proves by induction on $n$ that the $n$-th
derivative
\[
(k_\epsilon \ast f)^{(n)}(x)
\]
is a finite sum of terms of the form
\[
\frac{1}{\sqrt{2\pi}\epsilon}\int_{\mathbb R} \exp \big(
-{\frac{(u-x)^2}{2\epsilon^2}}\big)p(x-u)f(u)du,
\]
where $p$ is a polynomial. All limits,
\[
\lim_{|x|\rightarrow\infty}x^m (k_\epsilon \ast f)^{(n)}(x)=0
\]
are then shown, using the dominated convergence theorem, and all
the functions
\[
(k_\epsilon \ast f)(x)=\frac{1}{\sqrt{2\pi}\epsilon}\int_{\mathbb
R} \exp (-{\frac{(u-x)^2}{\epsilon^2}} )f(u)du
\]
are in the Schwartz space.
\mbox{}\qed\mbox{}\\

We now compute the adjoint operator $Q_\sigma^*$. Note that it is an
operator from $\mathbf L_2(\mathbb R,dx)$  into $\mathcal
S^\prime$, and therefore lies outside ${\mathbf L}_2(\mathbb
R,dx)$. We begin with a notation and a preliminary computation.
For $\phi\in\mathbf L_2(\mathbb R,dx)$ set
\[
c_n(\phi)=\int_{\mathbb R}\phi(y)h_n(y)dy.
\]
Then $(c_n(\phi))_{n\in\mathbb N}$ is in $\ell_2$. Indeed
\[
\|(c_n(\phi))_{n\in\mathbb N}\|_{\ell_2}^2=\|\psi\|^2_{\mathbf
L_2(\mathbb R,dx)}.
\]
We introduce the operator from $\mathbf L_2(\mathbb R,dx)$ into
$\mathbf L_2(d\sigma)$:
\[
(T_\sigma\phi)(u)=\sum_{n=0}^\infty c_n(\phi)e^{i\lambda_n u}.
\]
Clearly
\[
\|\phi\|_{\mathbf L_2(\mathbb R,dx)}=\|T_\sigma\phi\|_{\mathbf
L_2(d\sigma)}.
\]
\begin{Tm}
Let $\psi\in\mathcal S$ and $\phi\in\mathbf L_2(\mathbb R,dx)$.
Then,
\begin{eqnarray}
\label{T_h}
\langle Q_\sigma\psi,\phi\rangle_{\mathbf L_2(\mathbb
R,dx)}&=&\int_{\rm
sup~\sigma}\widehat{\psi}(u)\overline{(T_\sigma(\phi))}d\sigma(u)\\
\label{T_h2} &=&\int_{\mathbb R}\psi(y){\mathbf X(\phi)(y)})dy,
\end{eqnarray}
where
\begin{equation}
\label{Q*2}
(\mathbf X(\phi))(y)=\sum_{n=0}^\infty \langle
h_n,\phi\rangle_{\mathbf L_2(\mathbb
R,dx)}\widehat{\sigma}(y-\lambda_n).
\end{equation}
\end{Tm}

{\bf Proof:} We first prove \eqref{T_h}. In view of the formula
\eqref{Qpsi} for $Q_\sigma$, we have
\[
\begin{split}
\langle Q_\sigma\psi,\phi\rangle_{\mathbf L_2(\mathbb
R,dx)}&=\sum_{n=0}^\infty\left(\int_{\mathbb
R}\widehat{\sigma}(y-\lambda_n)\psi(y)dy\right)\overline{\left(\int_{\mathbb
R} \phi(x)h_n(x)dx\right)}\\
&\hspace{-2cm}=\sum_{n=0}^\infty\left(\int_{\mathbb
R}\left(\int_{\mathbb R}e^{iu(y-\lambda_n)}d\sigma(u)\right)
\psi(y)dy\right)\overline{\left(\int_{\mathbb R}
\phi(x)h_n(x)dx\right)}\\
&\hspace{-2cm}=\sum_{n=0}^\infty\left(\int_{\mathbb R}
\left(\int_{\mathbb R}e^{iuy}\psi(y)dy\right)
e^{-iu\lambda_n}d\sigma(u) \right)\overline{\left(\int_{\mathbb R}
\phi(x)h_n(x)dx\right)}\\
&\hspace{-2cm}=\int_{\mathbb R
}\widehat{\psi}(u)\left(\sum_{n=0}^\infty e^{-i\lambda_n
u}\overline{\left(\int_{\mathbb
R} \phi(x)h_n(x)dx\right)}\right)d\sigma(u)\\
&\hspace{-2cm}=\int_{\mathbb R}\widehat{\psi}(u)\cdot
\overline{T_\sigma(\phi)}d\sigma(u),
\end{split}
\]
where we have used Fubini's theorem for the third equality,
and the continuity of the inner product for the fourth equality.\\

We now turn to the second formula. The sequence $( \langle
h_n,\phi\rangle_{\mathbf L_2(\mathbb R,dx)})_{n\in\mathbb N_0}$
is in $\ell_2$. In view of \eqref{eq:equiv1}, the Cauchy-Schwarz
inequality implies that $\mathbf X(\phi)$ in \eqref{Q*2} converges
pointwise for every real $y$. We note that, in general, $\mathbf
X(\phi)\not\in\mathbf L_2(\mathbb R,dx)$. For $\phi$ and $\psi$
as in \eqref{T_h} we have:
\[
\begin{split}
\langle Q_\sigma\psi,\phi\rangle_{\mathbf L_2(\mathbb
R,dx)}&=\int_{\mathbb R}\left(\sum_{n=0}^\infty\left(\int_{\mathbb
R}\widehat{\sigma}(y-\lambda_n)\psi(y)dy\right)h_n(x)\right)
\overline{\phi(x)}dx\\
&=\sum_{n=0}^\infty\left(\int_{\mathbb
R}\widehat{\sigma}(y-\lambda_n)\psi(y)dy\right)\left(\int_{\mathbb
R}h_n(x)\overline{\phi(x)}dx\right)\\
&=\int_{\mathbb
R}\psi(y)\left(\sum_{n=0}^\infty\widehat{\sigma}(y-\lambda_n)
\overline{\int_{\mathbb R} \phi(x)h_n(x)dx}\right)dy.
\end{split}
\]
To obtain the second equality, we note the following: Write
\[
\int_{\mathbb R}\widehat{\sigma}(y-\lambda_n)\psi(y)dy=
\int_{\mathbb
R}\frac{\widehat{\sigma}(y-\lambda_n)}{y^2+1}((y^2+1)\psi(y))dy.
\]
Using the Cauchy-Schwarz inequality, we see that
\[
\big| \int_{\mathbb
R}\widehat{\sigma}(y-\lambda_n)\psi(y)dy\big|^2\le\left(\int_{\mathbb
R} \frac{|\widehat{\sigma}(y-\lambda_n)|^2}{(y^2+1)^2}\right)
\left(\int_{\mathbb R}(y^2+1)^2)|\psi(y)|^2dy\right).
\]
In view of \eqref{eq:equiv1},
\[
\sum_{n=0}^\infty\left(\int_{\mathbb
R}\widehat{\sigma}(y-\lambda_n)\psi(y)dy\right)h_n(x)
\]
belongs to $\mathbf L_2(\mathbb R,dx)$, and we use the continuity
of the scalar product. Furthermore we have used the dominated
convergence theorem to obtain the third equality.
\mbox{}\qed\mbox{}\\

Equations \eqref{T_h} and \eqref{T_h2} allow to compute
$Q_\sigma^*$:

\begin{Tm}
The domain of $Q_\sigma^*$ is the Lebesgue space $\mathbf L_2(\mathbb R,
dx)$. For $\phi\in\mathbf L_2(\mathbb R, dx)$, $Q_\sigma^*\phi$ is the
tempered distribution defined through
\[
\langle Q_\sigma^*(\phi),\psi\rangle_{\mathcal S^\prime,\mathcal
S}= \int_{\rm
sup~\sigma}\overline{\widehat{\psi}(u)}(T_\sigma(\phi))d\sigma(u).
\]
Equivalently, $Q_\sigma^*(\phi)$ is the tempered distribution
defined by the function $\mathbf X(\phi)$, that is, with some abuse of
notation
\[
Q_\sigma^*(\psi)(y)= \sum_{n=0}^\infty \langle
h_n,\phi\rangle_{\mathbf L_2(\mathbb
R,dx)}\widehat{\sigma}(y-\lambda_n).
\]
\label{tm56}
\end{Tm}

The second representation for $Q_\sigma^*$ has an important
consequence:
\begin{Tm}
It holds that
\[
\ker Q_\sigma^*=\left\{0\right\}.
\]
\end{Tm}

It follows that the operator $Q_\sigma^*Q_\sigma$ is a continuous
and bounded operator from $\mathcal S$ into $\mathcal S^\prime$.
We now provide two formulas for this operator.\\

\begin{Tm}
\label{tm58}
Let $\psi$ and $\phi$ be in $\mathcal S$. Then,
\begin{equation*}
\begin{split}
\hspace{-1cm}
\langle (Q_\sigma^*Q_\sigma)\phi,\psi
\rangle_{\mathcal S,\mathcal S^\prime}&=\\
&\hspace{-2cm} =\sum_{n=0}^\infty \left(\int_{\mathbb R
}\widehat{\phi}(u)e^{-i\lambda_n
u}d\sigma(u)\right)\overline{\left(\int_{\mathbb R}\widehat{\psi}(u)
e^{-i\lambda_n u}d\sigma(u)\right)}\\
&\hspace{-2cm} =\int_{\mathbb R}
\widehat{\phi}(u)\overline{\widehat{\psi}(u)}d\sigma(u).
\end{split}
\end{equation*}
\end{Tm}

{\bf Proof:} By definitions of $Q_\sigma\phi$ and of $Q_\sigma^*$ we have
\[
\begin{split}
\langle (Q_\sigma^*Q_\sigma)(\phi),\psi\rangle_{\mathcal
S,\mathcal S^\prime}&=\int_{\mathbb R}\overline{\widehat{\psi}(u)}
\left(T_\sigma(Q_\sigma\phi)\right)\\
&\hspace{-2cm}= \sum_{n=0}^\infty \int_{\mathbb
R}\overline{\widehat{\psi}(u)}\left( \int_{\mathbb
R}\widehat{\sigma}(y-\lambda_n)\phi(y)dy\right)e^{i\lambda_n
u}d\sigma(u)\\
&\hspace{-2cm}=\sum_{n=0}^\infty
\left(\int_{\mathbb R}\overline{\widehat{\psi}(u)e^{-i\lambda_n
u}}d\sigma(u)\right)\left(
\int_{\mathbb R}\widehat{\sigma}(y-\lambda_n)\phi(y)dy\right).
\end{split}
\]
But
\[
\begin{split}
\int_{\mathbb
R}\widehat{\sigma}(y-\lambda_n)\phi(y)dy&=\int_{\mathbb R}\left(
\int_{\mathbb
R}e^{i(y-\lambda_n)v}d\sigma(v)\right)\phi(y)dy\\
&=\int_{\mathbb R}\widehat{\phi}(y)e^{-i\lambda_n v}d\sigma(v)
\end{split}
\]
where we have used Fubini's theorem. This concludes the proof.
\mbox{}\qed\mbox{}\\

\begin{Rk}
{\rm While the operator $Q_\sigma$ (see Theorems \ref{tm56}
through  \ref{tm58}) is well defined as an unbounded linear
operator in the Hilbert space $\mathbf L_2(\mathbb R,dx)$, the
other two operators $Q_\sigma^*$ and $Q_\sigma^*Q_\sigma$ are not.
The reason is that the range of $Q_\sigma^*$ is not contained in
$\mathbf L_2(\mathbb R,dx)$. In fact,
\[
(Q_\sigma^*h_n)(x)=\widehat{\sigma}(x-\lambda_n),
\]
where $\widehat{\sigma}$  is the infinite product expression
\eqref{sigmarho}. It can be shown that $x\mapsto
\widehat{\sigma}(x)$ is not in $\mathbf L_2(\mathbb R,dx)$; so as
an $\mathbf L_2(\mathbb R,dx)$-operator, $Q_\sigma$ is not
closable (its adjoint, computed in $\mathbf L_2(\mathbb R, dx)$,
does not have dense domain! Nonetheless $Q_\sigma^*$ in the
extended sense maps $\mathbf L_2(\mathbb R, dx)$ into $\mathcal
S^\prime$). The use of the ambient space $\mathcal S^\prime$ of
tempered distributions is essential. We illustrate the above
discussion with the following diagrams:

\[\renewcommand{\arraystretch}{1.5}
\begin{array}{ccccc}
\mathcal S&\stackrel{i}{\longhookrightarrow}
&\mathbf L_2(\mathbb R,dx)& \stackrel{i^*}{\longhookrightarrow}
&\mathcal S^\prime\\
&\stackrel{Q_\sigma}{\searrow}& &\stackrel{Q_\sigma^*}{\nearrow}\\
\mathcal S&\stackrel{i}{\longhookrightarrow}  &\mathbf L_2(\mathbb
R,dx)&\stackrel{i^*}{\longhookrightarrow}&\mathcal S^\prime
\end{array}.\]

In the following diagram, ${\rm dom}~Q_\sigma^*$ is only a small
subspace of $\mathbf L_2(\mathbb R,dx)$:

\[\renewcommand{\arraystretch}{1}
\begin{array}{ccccc}
& &\mathbf L_2(\mathbb R,dx)& \stackrel{i^*}{\longhookrightarrow}
&\mathcal S^\prime\\
&{\text{\tiny restricted}}\stackrel{Q_\sigma^*}{\nearrow}
& &\stackrel{Q_\sigma^*}{\nearrow}{\text{\tiny unbounded}}\\
{\rm
Dom}~Q_\sigma^*\hspace{-.7cm}&\hspace{5mm}\stackrel{i}{\longhookrightarrow}
&\mathbf L_2(\mathbb R,dx)&
\end{array}.\]

}
\end{Rk}
\section{The processes $X_\sigma$ and $W_\sigma$}
\label{secXW}
\setcounter{equation}{0}
In this section we build the Gaussian process with covariance
function $K_\sigma$. First recall that, thanks to Proposition
\ref{10t}, the domain of the operator $Q_\sigma$ has been extended
to include the functions $1_{[0,t]}$. We begin with:

\begin{Tm}
Let $Q_\sigma$ be as \eqref{secQ}. Then, for every $t,s\in\mathbb
R$,
\[
\langle Q_\sigma 1_{[0,t]}, Q_\sigma 1_{[0,s]}\rangle_{{\mathbf
L_2(\mathbb R,dx)}}= \int_{\mathbb R
}\frac{e^{iut}-1}{u}\frac{e^{-ius}-1}{u}d\sigma(u).
\]
\end{Tm}

{\bf Proof:} We take $(s_n)_{n\in\mathbb N}$ and
$(t_n)_{n\in\mathbb N}$ two sequences of elements of $\mathcal S$
with Fourier transforms $(\widehat{s_n})_{n\in\mathbb N}$ and
$(\widehat{t_n})_{n\in\mathbb N}$ converging in the supremum norm
respectively to $\chi_t$ and $\chi_s^*$. Then $Q1_{[0,t]}$ and
$Q1_{[0,s]}$ are the limit in $\mathbf L_2(\mathbb R,dx)$ of the
sequences $(Qs_n)_{n\in\mathbb N}$ and $(Qt_n)_{n\in\mathbb N}$
respectively. Hence,
\[
\begin{split}
\langle Q_\sigma 1_{[0,t]}, Q_\sigma 1_{[0,s]}\rangle_{{\mathbf
L_2(\mathbb R,dx)}} &= \lim_{n,m\rightarrow\infty}\langle
Q_\sigma s_n,
Q_\sigma t_m\rangle_{{\mathbf L_2(\mathbb R,dx)}}\\
&=\lim_{n,m\rightarrow\infty}\int_{\mathbb R}
\widehat{s_n}(u)\widehat{t_n}d\sigma(u)\\
&= \int_{\mathbb R
}\frac{e^{iut}-1}{u}\frac{e^{-ius}-1}{u}d\sigma(u).
\end{split}
\]
\mbox{}\qed\mbox{}\\

Recall that we have denoted by $\widetilde{f}$ the natural
isometric imbedding \eqref{ftilde} of $\mathbf L_2(\mathbb R,dx)$
into the white noise space. We set
\begin{equation}
\label{ZN}
Z_n=\widetilde{h_n},\quad n=0,1,\ldots
\end{equation}
The $Z_n$ are independent, identically distributed
$\mathcal N(0,1)$ random variables, and represented in white noise space.\\

We arrive at the following decomposition:
\begin{equation}
\label{QXt}
\begin{split}
X_\sigma(t)&=\widetilde{Q_\sigma(1_{[0,t]})}\\
&=\sum_{n=0}^\infty\left(\int_0^t\widehat{\sigma}(y-\lambda_n)dy\right)Z_n.
\end{split}
\end{equation}
\begin{Tm}
\label{tm:kl}
It holds that
\begin{equation}
\label{Xsigmat-s} \|X_\sigma(t)-X_\sigma(s)\|_{\mathcal W}\le
|t-s|,\quad t,s\in\mathbb R.
\end{equation}
The function $t\mapsto X_\sigma(t)$ is differentiable in $\mathcal
W$ (white noise space), and its derivative is given by
\begin{equation}
\label{eq:Wsigma}
W_\sigma(t)=\sum_{n=0}^\infty\widehat{\sigma}(t-\lambda_n)Z_n
\end{equation}
\end{Tm}

\begin{Rk}
{\rm Formulas \eqref{QXt} and \eqref{eq:Wsigma} are analogous to,
but different from Karhunen-Lo\`eve expansions; see e.g.,
\cite{JoSo09a}.}
\end{Rk}

{\bf Proof of Theorem \ref{tm:kl}:} We first prove
\eqref{Xsigmat-s}. We have
\[
\begin{split}
E_\sigma[|X_\sigma(t)-X_\sigma(s)|^2]&=
\|X_\sigma(t)-X_\sigma(s)\|_{\mathcal W}^2\\
&=\sum_{n=0}^\infty
\big|\int_s^t\widehat{\sigma}(y-\lambda_n)dy\big|^2\\
&\le \sum_{n=0}^\infty(\int_s^t 1^2
dy)(\int_s^t|\widehat{\sigma}(y-\lambda_n)|^2dy)\\
&=(t-s)\int_s^t(\sum_{n=0}^\infty|\widehat{\sigma}(y-\lambda_n)|^2)dy\\
&=(t-s)^2,
\end{split}
\]
where we have used the Cauchy-Schwarz inequality and
\eqref{eq:equiv1} for the second and fourth equalities,
respectively.\\

We remark that, in view of \eqref{eq:equiv1},
$W_\sigma(t)\in\mathcal W$ for every real $t$. We have
\[
\begin{split}
\frac{X_\sigma(t)-X_\sigma(s)}{t-s}-W_\sigma(t)&=\frac{\sum_{n=0}^\infty
\int_s^t\widehat{\sigma}(y-\lambda_n)dy Z_n}{t-s}-\\
&\hspace{10mm}-
\sum_{n=0}^\infty \widehat{\sigma}(t-\lambda_n)dy Z_n\\
&=\frac{\sum_{n=0}^\infty
\left\{\int_s^t\left(\widehat{\sigma}(y-\lambda_n)-
\widehat{\sigma}(t-\lambda_n)\right)dy\right\}
Z_n}{t-s}\\
&=\frac{\sum_{n=0}^\infty \left\{\int_s^t\langle
e_y-e_t,e_{\lambda_n}\rangle_{\mathbf L_2(d\sigma)} dy\right\}
Z_n}{t-s}.
\end{split}
\]
Hence, using Parseval's equality and \eqref{ineq:ets}, we
obtain
\[
\begin{split}
\left\|
\frac{X_\sigma(t)-X_\sigma(s)}{t-s}-W_\sigma(t)\right\|^2_{\mathcal
W}&=\frac{\sum_{n=0}^\infty\big|\int_s^t \langle
e_y-e_t,e_{\lambda_n}\rangle_{\mathbf
L_2(d\sigma)}dy\big|^2}{(t-s)^2}\\
&\le\sum_{n=0}^\infty\frac{1}{t-s}\int_s^t |\langle
e_y-e_t,e_{\lambda_n}\rangle_{\mathbf L_2(d\sigma)}|^2dy\\
&=\frac{1}{t-s}\int_s^t\|e_y-e_t\|^2_{\mathbf L_2(d\sigma)}dy\\
&\le \frac{1}{t-s}\int_s^t(y-t)^2dy\\
&\le \frac{(t-s)^2}{3}.
\end{split}
\]
\mbox{}\qed\mbox{}\\

\begin{Tm}
The derivative process $W_\sigma$ is continuous in the
$\|\cdot\|_{\mathcal W}$ norm. It is furthermore stationary and
of constant variance,
\begin{equation}
\label{conscv}
E[|W_\sigma(t)|^2]\equiv 1.
\end{equation}
\end{Tm}

{\bf Proof:}
\[
\begin{split}
\|W_\sigma(t)-W_\sigma(s)\|_{\mathcal
W}^2&=\sum_{n=0}^\infty|\widehat{\sigma}(t-\lambda_n)-
\widehat{\sigma}(s-\lambda_n)|^2\\
&= \sum_{n=0}^\infty|\langle e_t-e_s,
e_{\lambda_n}\rangle_{\mathbf L_2(d\sigma)}|^2\\
&=\|e_t-e_s\|_{\mathbf L_2(d\sigma)}^2\\
&=2(1-\langle e_t,e_s\rangle_{\mathbf L_2(d\sigma)})\\
&=2(1-\widehat{\sigma}(t-s))\\
&=2\left(1-\prod_{n=1}^\infty
\cos\left(\frac{t-s}{4^n}\right)\right).
\end{split}
\]
Continuity and stationarity follow from the above chain of
inequalities, together with the fact that $\left\{W_\sigma(t)\right\}$
is a Gaussian process. Equation \eqref{conscv} follows from
\eqref{eq:equiv1}.
\mbox{}\qed\mbox{}\\

We now turn to another type of representation for $X_\sigma$:

\begin{Tm}
Let
\[
H(\w,u)=\sum_{n=0}^\infty \widetilde{h_n}(\w)e^{i\lambda_n
u}\in\mathcal W\otimes\mathbf L_2(d\sigma).
\]
Then
\[
(X_\sigma(t))(\w)=\int_{\mathbb R
}H(\w,u)\frac{e^{-iut}-1}{u}d\sigma(u)
\]
\end{Tm}

{\bf Proof:} Using Fubini's theorem we have
\[
\begin{split}
\int_{\mathbb R}e^{i\lambda_nu}\frac{e^{-iut}-1}{u}d\sigma(u)&=
\int_{\mathbb R}e^{i\lambda_nu}(\int_0^te^{-iuv}dv)d\sigma(u)\\
&=\int_0^t\int_{\mathbb R}e^{i(\lambda_n-v)u}d\sigma(u)\\
&=\int_0^t\widehat{\sigma}(\lambda_n-v).
\end{split}
\]
This together with \eqref{QXt} leads to the required conclusion.
\mbox{}\qed\mbox{}\\

\section{The Wick-Ito integral}
\setcounter{equation}{0}
\label{wickito}
In this section we establish a stochastic integration formula for
the general class of stationary increment processes considered
here. An extension of Ito's formula is addressed in section
\ref{itoform}. More precisely, see Theorem \ref{bastille}, with
the use of {\it a priori} estimates and of a Wick calculus (from
weighted symmetric Fock spaces), we extend and sharpen earlier
computations of Ito-stochastic integration developed originally
only for the special case of stationary increment processes
having absolutely continuous spectral measures. We further obtain
in the subsequent section an associated Ito formula (Theorem 8.2).\\

The main result of this section is the following theorem, which
is the counterpart of \cite[Theorem 5.1]{ptrf}. The fact that the
derivative process $W_\sigma(t)$ is $\mathcal W$-valued (rather than lying in the
larger Kondratiev space) allows to
get a sharper statement. The proof follows the strategy of \cite{ptrf}
and hence is only outlined.

\begin{Tm}
\label{bastille}
Let $Y(t)$, $t\in[a,b]$ be an $S_{-1}$-valued
function, continuous in the strong topology of $S_{-1}$. Then,
there exists a $p\in\mathbb N$ such that the function $t\mapsto
Y(t)\lozenge W_\sigma(t)$ is $\mathcal H_p^\prime$-valued, and
\begin{equation*}
\int_a^b Y(t,\omega)\lozenge W_\sigma(t)dt
=\lim_{\left|\Delta\right|\to 0} \sum_{k=0}^{n-1}
Y(t_k,\omega)\lozenge
\left(X_\sigma(t_{k+1})-X_\sigma(t_k)\right),
\end{equation*}
where the limit is in the $\mathcal H_p^\prime$ norm, with
$\Delta: a=t_0<t_1<\cdots <t_n=b$ a partition of the interval
$[a,b]$ and $\left|\Delta\right|=\max_{0\leq k \leq
n-1}(t_{k+1}-t_k)$.
\end{Tm}

{\bf Proof:} We proceed in a number of steps.\\

STEP 1:  {\it There exists a $p\in\N$ such that
$Y(t)\in\mathcal{H}_{p}^\prime$ for all $t\in [a,b]$, being
uniformly continuous from $[a,b]$ into
$\mathcal{H}_{p}^\prime$.}\\

This is proved in \cite[STEP 2 of the Proof of Theorem 5.1]{ptrf}.\\





STEP 2: {\it The function $t\mapsto Y(t)\lozenge W_\sigma(t)$ is
continuous over $[a,b]$}.\\

Here and in the sequel, we set $\|\cdot\|_p\stackrel{\rm def.}{=}
\|\cdot\|_{{\mathcal H}_p^\prime}$ to simplify notation. Using
V\r{a}ge's inequality \eqref{vage}, it follows that, for $p>1$,
\[
\begin{split}
\left\|Y(t)\lozenge W_\sigma(t)
-Y(s)\lozenge W_\sigma(s)\right\|_p&\le\\
&\hspace{-4cm}\le \left\|(Y(t)-Y(s))\lozenge
W_\sigma(t)\right\|_p+
\left\|Y(s)\lozenge (W_\sigma(t)-W_\sigma(s))\right\|_p\\
&\hspace{-4cm}\le A(p)\|Y(t)-Y(s)\|_p\|W_\sigma(s)\|_{0}+\\
&\hspace{-3.5cm}+A(p)\|Y(s)\|_p\|W_\sigma(t)-W_\sigma(s)\|_{0}
\end{split}
\]
where $A(p)$ is defined by \eqref{vage111}. We conclude the proof
of STEP 2 by observing that
\[
\|\cdot\|_{\mathcal H_0}\le \|\cdot\|_{\mathcal W},
\]
which implies the continuity of the function $Y(t)\lozenge W_\sigma(t)$ in the
norm $\|\cdot\|_{\mathcal W}$.\\

In view of Step 2, the integral $\int_a^b Y(t)\lozenge
W_\sigma(t)dt$ makes sense as a Riemann integral of a
continuous Hilbert space valued function.\\

STEP 3: {\it Let $\Delta$ be a partition of the interval $[a,b]$.
We now compute an estimate for}
\[
\begin{split}
\int_a^b Y(t)\lozenge W_\sigma(t)dt-\sum_{k=0}^{n-1}
Y(t_k)\lozenge
\left(X_\sigma(t_{k+1})-X_\sigma(t_k)\right)&=\\
&\hspace{-6cm}= \sum_{k=0}^{n-1}
\left(\int_{t_k}^{t_{k+1}}(Y(t)-Y(t_k)) \lozenge
W_\sigma(t)dt\right).
\end{split}
\]

As for steps 1 and 2, we closely follow \cite{ptrf}. Let $p$ be as
in Step 2, and set $\epsilon>0$. Since $Y$ is uniformly
continuous on $[a,b]$, there exists an $\eta>0$ such that
\[
|t-s|<\eta\Longrightarrow \|Y(t)-Y(s)\|_p<\epsilon.
\]

Set
\[
\tilde{C}=\max_{s\in[a,b]}\|W_\sigma(s)\|_{0}\quad{\rm and}\quad
A=A(p-N-3).
\]

Let $\Delta$ be a partition of $[a,b]$ with
\[
|\Delta|=\max\left\{|t_{k+1}-t_k|\right\}<\eta.
\]

We then have:

\[
\begin{split}
\left\|\sum_{k=0}^{n-1} \left(\int_{t_k}^{t_{k+1}}(Y(t)-Y(t_k))
\lozenge W_\sigma(t)dt \right)\right\|_p\leq\\
&\hspace{-7cm}\le \sum_{k=0}^{n-1}
\left(\int_{t_k}^{t_{k+1}}\left\|(Y(t)-Y(t_k)) \lozenge
W_\sigma(t)\right\|_p dt\right)\\
&\hspace{-7cm}\leq A\sum_{k=0}^{n-1}
\left(\int_{t_k}^{t_{k+1}}\left\|(Y(t)-Y(t_k))\right\|_p\left\|
W_\sigma(t)\right\|_{0} dt\right)\\
&\hspace{-7cm}\leq \tilde{C}A\sum_{k=0}^{n-1}\int_{t_k}^{t_{k+1}}
\left\|(Y(t)-Y(t_k))\right\|_pdt\\
&\hspace{-7cm}\le \epsilon\tilde{C}A(b-a),
\end{split}
\]
which completes the proof of Step 3 and the proof of the Theorem.
\mbox{}\qed\mbox{}\\


\section{An Ito formula}
\setcounter{equation}{0}
\label{itoform}
We extend the classical Ito's formula to the present setting. Our
present wider context for these stochastic processes entails
important analytical points: Singular measures of fractal
dimension, and singular operators, extending beyond the Hilbert
space $\mathbf L_2(\mathbb R,dx)$. This in turn brings to light
new aspects of Ito calculus which we detail below.
In addition to the examples in Section \ref{sec:self} (affine IFS
measures), we further offer two examples in sect \ref{two_ex}
below: the periodic Brownian bridge, and the Orenstein-Uhlenbeck
processes.

\begin{La}
The function
\[
r(t)=\|Q_\sigma1_{[0,t]}\|_{{\mathbf L_2(\mathbb R)}}^2
\]
is absolutely continuous with respect to the Lebesgue measure.
\end{La}

{\bf Proof:} By \eqref{normQ} we have
\[
\begin{split}
\|Q_\sigma1_{[0,t]}\|_{{\mathbf L_2(\mathbb R)}}^2 &=\int_{\mathbb R}
|\chi_t(u)|^2d\sigma(u)\\
&=2\int_{\mathbb R}\frac{1-\cos(tu)}{u^2}d\sigma(u).
\end{split}
\]
Since the support of $d\sigma$ is bounded, the dominated
convergence theorem allows to show that $r(t)$ is differentiable
and that its derivative is given by
\[
r^\prime(t)=2\int_{\mathbb R}\frac{\sin(tu)}{u}d\sigma(u).
\]
\mbox{}\qed\mbox{}\\

\begin{Tm}
Let $f:\R\to\R$ be a $C^2(\mathbb R)$ function. Then
\begin{equation}
\label{oberkampf}
\begin{split}
f(X_\sigma(t))&=f(X_\sigma (t_0 ))+\int_{t_0}^t f'(X_\sigma(s))
\lozenge W_\sigma(s)ds+\\
&\hspace{5mm}+
\frac{1}{2}\int_{t_0}^tf^{\prime\prime}(X_\sigma(s)) r^\prime(s)
ds,\quad t_0 <t \in\R,
\end{split}
\end{equation}
where the equality is in the $P$-almost sure sense.
\end{Tm}

{\bf Proof:}  We prove for $t>t_0 =0$. The proof for any other
 interval in $\R$ is essentially the same.
We divide the proof into a number of steps. Step 1-Step 5 are
constructed so as to show that \eqref{oberkampf} holds, $\forall
t >0$, for  Schwartz functions. This enables the extension to
$C^2$ functions with compact support, with the
equality holding in the ${\mathcal H}^\prime_p$ sense. This
implies its validity in the $P$-a.s. sense (actually, holding
$\forall \omega\in\Omega$),  hence, setting the ground for the
concluding step, in which the result is extended
to hold for all $C^2$  functions $f$.\\

STEP 1: {\it For every $(u,t)\in\mathbb R^2$, it holds that
\[
e^{iuX_\sigma(t)}\in\mathcal W,
\]
and
\begin{equation}
e^{iuX_\sigma(t)}\lozenge W_\sigma(t)\in\mathcal H^\prime_{2}.
\label{leipzig2010}
\end{equation}

}

Indeed, since $X_\sigma$ is real, we have
\[
|e^{iuX_\sigma(t)}|\le 1,\quad \forall u,t\in\mathbb R,
\]
and hence $e^{iuX_\sigma(t)}\in\W$. Since $\W\subset \mathcal
H_{2}^\prime$, we have in particular that $W_\sigma(t)\in \mathcal
H_{2}^\prime$ for all $t\in\mathbb R$, it follows from V\r{a}ge's
inequality
\eqref{vage} that \eqref{leipzig2010} holds.\\

The aim of the following two steps is formula \eqref{oberkampf}
for exponential functions. For $\alpha\in\R$ we set:
\[
g(x)=\exp(i\alpha x).
\]
The proofs are as in \cite{ptrf}, taking into account that $r(t)$
is absolutely continuous with respect to Lebesgue measure and are
omitted.\\

STEP 2: {\it It holds that
\begin{equation}
\label{gprimem}
 g'(X_\sigma(t))=i\alpha g(X_\sigma(t))\lozenge
W_\sigma(t)+\frac{1}{2} (i\alpha)^2g(X_\sigma(t))r^\prime (t).
\end{equation}}

STEP 3: {\it Equation \eqref{oberkampf} holds for exponentials.}\\

In the following two steps, we prove \eqref{oberkampf} to hold
for Schwartz functions.\\

STEP 4: {\it The function $(u,t)\mapsto e^{iuX_\sigma(t)}\lozenge
W_\sigma(t)$ is continuous from $\mathbb R^2$ into $\mathcal
H^\prime_{2}$.}\\

We first recall that the norm in $\mathcal H_p^\prime$ is denoted
by $\|\cdot\|_{p}$. See \eqref{michelle}. The particular case $p=0$ in
\eqref{michelle} gives in particular
\[
\|f\|_{0}^2=\sum_{\alpha\in\ell}|f_\alpha|^2.
\]
The structure of $\mathcal H^\prime_0$ has been studied in
\cite[Section 7]{alpay-2008}.\\

Recall now that the function $t\mapsto X_\sigma(t)$ is
continuous, and even uniformly continuous, from $\mathbb R$ into
$\mathcal W$, and hence from $\mathbb R$ into $\mathcal
H_{p}^\prime$ for any $p\in\mathbb N_0$ since
\[
\|u\|_{p}
\le\|u\|_{\W},\quad{\rm for}\quad
u\in\W.
\]
The function $(u,t)\mapsto e^{iuX_\sigma(t)}$ is in particular
continuous from $\mathbb R^2$ into $\mathcal H^\prime_{2}$.
Furthermore, using V\r{a}ge's inequality \eqref{vage} we have:
\[
\begin{split}
\|e^{iu_1X_\sigma(t_1)}\lozenge
W_\sigma(t_1)-e^{iu_2X_\sigma(t_2)}\lozenge
W_\sigma(t_2)\|_{\mathcal H^\prime_{2}}&\le\\
&\hspace{-25mm}\le \|
(e^{iu_1X_\sigma(t_1)}-e^{iu_2X_\sigma(t_2)})
\lozenge W_\sigma(t_1)\|_{\mathcal H^\prime_{2}}+\\
&\hspace{-20mm}+
 \| e^{iu_1X_\sigma(t_1)}\lozenge(W_\sigma(t_2)-
W_\sigma(t_1))\|_{\mathcal H^\prime_{2}}\\
&\hspace{-65mm}\le A(2)\| (e^{iu_1X_\sigma(t_1)}
-e^{iu_2X_\sigma(t_2)})
\|_{\mathcal H^\prime_{2}}\cdot\| W_\sigma(t_1)\|_{\mathcal H^\prime_{0}}+\\
&\hspace{-60mm}+ A(2) \| e^{iu_1X_\sigma(t_1)}\|_{\mathcal
H^\prime_{2}}\cdot\|(W_\sigma(t_2)- W_\sigma(t_1))\|_{\mathcal
H^\prime_{0}},
\end{split}
\]
where $A(2)$ is defined by \eqref{vage111}. This completes the
proof of STEP 4 since
\[
\|\cdot\|_2\le\|\cdot\|_0\le\|\cdot\|_{\mathcal W}
\]
and in particular $t\mapsto W_\sigma(t)$ is continuous in the norm
of $\mathcal H_{0}^\prime$ and $(u,t)\mapsto e^{iuX_\sigma(t)}$ is
continuous in the norm of
$\mathcal H_{2}^\prime$.\\

STEP 5: {\it  \eqref{oberkampf} holds for $f$ in the Schwartz
space.}\\

The reminder of the proof is exactly as Steps 6-9 in the corresponding
proof of \cite[Theorem 6.1]{ptrf}, and hence omitted.
\mbox{}\qed\mbox{}\\

\section{Two examples}
\label{two_ex}
\setcounter{equation}{0}
In this section, to contrast the measures considered above, we now
consider two example where the measure $ \sigma$ has an unbounded
support.\\

{\bf Example (The periodic Brownian bridge)}.  Take
\begin{equation}
\label{browbri} \sigma(u)=\sum_{n=0}^\infty \delta(u-2n),
\end{equation}
that is the measure with support on the even positive integers
with mass equal to $1$ at each of these points. Then, it follows
from the formula \eqref{rreal} on $r(t)$ that
\[
r(t)=\pi\sum_{n=1}^\infty \frac{1-\cos(2nt)}{(2n)^2}.
\]
Note that,
\[
t(\pi-t)=\pi\sum_{n=1}^\infty \frac{1-\cos(2nt)}{(2n)^2},\quad
t\in[0,2\pi].
\]
In view of the preceding equality, we call the associated process
$\left\{X_\sigma(t)\right\}$ the {\it periodic Brownian bridge
over $[0,\pi]$}. We have
\[
X_\sigma(t)=\sqrt{\frac{\pi}{2}}\sum_{n=1}^\infty\frac{\sin
(nt)}{n}Z_n,
\]
where, as in  \eqref{ZN}, $Z_n=\widetilde{h_n}$.\\

We note that, by construction,  for every $t$, $X_\sigma(t)$
belongs to the white noise space. On the other hand, the power
series
\[
W_\sigma(t)=\sqrt{\frac{\pi}{2}}\sum_{n=1}^\infty\cos(nt)Z_n
\]
converges only in the Kondratiev space. More precisely,
\begin{Tm}
Let $\sigma$ be given by \eqref{browbri}. For every $t$, we have
that $W_\sigma(t)\in \mathcal H^\prime_2$ and in the topology of
$\mathcal H^\prime_4$
\begin{equation}
\label{eq:deriv} \lim_{s\rightarrow t}
\frac{X_\sigma(t)-X_\sigma(s)}{t-s}=W_\sigma(t).
\end{equation}
\end{Tm}

{\bf Proof:} The fact that $W_\sigma(t)$ belongs to $\mathcal
H^\prime_2$ follows from definition \eqref{michelle} since
\[
\|Z_n\|^2_2=(2n)^{-2},
\]
and since the $Z_n$ are mutually orthogonal in $\mathcal
H_2^\prime$. We now turn to \eqref{eq:deriv}. We have
\[
\begin{split}
\frac{X_\sigma(t)-X_\sigma(s)}{t-s}-W_\sigma(t)&=\sum_{n=0}^\infty
\left(\frac{\int_s^t(\cos(nu)-\cos(nt))du}{t-s}\right)Z_n.
\end{split}
\]
But
\[
\begin{split}
\Big|\frac{\int_s^t(\cos(nu)-\cos(nt))du}{t-s}\Big|&\le\Big|\frac{\int_s^t
(nu-ns)n\sin(nv)du}{t-s}\Big|,\quad\mbox{{\rm for some}}\,\,
v\in(s,t)\\
&\le n^2\Big|\frac{\int_s^t(u-s)du}{t-s}\Big|\\
&=\frac{n^2|t-s|}{2}
\end{split}
\]
and hence the limit goes to $0$ in the $\mathcal H^\prime_4$ norm
since
\[
\|Z_n\|^2_4=(2n)^{-4},
\]
and the $Z_n$ are mutually orthogonal in $\mathcal H_4^\prime$.
\mbox{}\qed\mbox{}\\

{\bf Example (The Ornstein-Uhlenbeck process)}. This is the
solution of ths stochastic differential equation
\[
dX(t)=\theta(\mu-X(t))dt+\alpha dB(t),\quad t\ge 0,
\]
where $B$ is a Brownian motion and $\theta\not=0, \mu$ and
$\sigma$ are parameters. We have
\[
X(t)=\mu+\frac{\alpha}{\sqrt{2\theta}}W(e^{2\theta t})e^{-\theta
t},
\]
and
\[
\begin{split}
E[X(t)]&=e^{-\theta t}+\mu(1-e^{-\mu t}),\\
E[(X(t)-E(X(t))(X(s)-E(X(s))]&=\frac{\alpha^2}{2\theta}
e^{-\theta(t+s)}(e^{2\theta s\wedge t}-1).
\end{split}
\]
\begin{Tm}
The centered Ornstein-Uhlenbeck process
$\left\{X(t)-E[X(t)]\right\}$ is a stationary increment Gaussian
process. Furthermore
\[
\begin{split}
E[|X(t)-E[X(t)]|^2]&=\frac{\alpha^2}{2\theta}(1-e^{-2\theta t})\\
&=\int_{\mathbb R} \frac{1-\cos(tu)}{u^2}d\sigma(u),
\end{split}
\]
with
\[
d\sigma(u)=\frac{\alpha^2}{2\pi\theta} \frac{\theta u^2
du}{\theta^2+u^2}.
\]
\end{Tm}



\section{The operator $Q_\sigma$ for more general spectral
measures $d\sigma$} \setcounter{equation}{0}
\label{QG}
The purpose of this section is to the contrast the difference
between a singular measure or not. This is a crucial distinction in passing
from the given measure  $\sigma$ to the associated process
$X_\sigma(t)$ ; i.e., in solving the inverse problem.  The two
cases are: $(i)$ $\sigma$  is assumed absolutely continuous with
respect to Lebesgue measure; versus: $(ii)$ $\sigma$ is singular.
This distinction results in a dichotomy for the induced operators
in $\mathbf L_2(\mathbb R, dx)$, so for the two operator
questions.
 In case $(i)$ we have a Radon-Nikodym derivative $m$  , and the induced operator is
 $T_m$, a selfadjoint convolution operator in the Hilbert space
 $\mathbf L_2(\mathbb R, dx)$  with the Schwartz space $\mathcal  S$ as dense
 domain; see \eqref{tm}.
 In the second case $(ii)$ it is not, even existence is subtle.
Now rather the induced operator is our operator  $Q_\sigma$  from
Section \ref{secQ}. But it is much more subtle, and below we
address some of the technical points omitted in Section
\ref{secQ}: $(a)$ There is now more than one choice for
$Q_\sigma$;
 and $(b)$ none of the choices will be closable operators,
 referring just to the Hilbert space $\mathbf L_2(\mathbb R, dx)$.
$(c)$ Nonetheless, by working within the environment of the
Gelfand triples (of Gaussian random fields), we are still able to
make precise the two operators  $Q_\sigma$ and the corresponding
adjoint  $Q_\sigma^*$. For the singular case, i.e.,  case $(ii)$,
the Gelfand triple  is thus essential in justifying our
construction of the process
$\left\{X_\sigma(t)\right\}$, existence and related properties.\\

One common point between the previous work \cite{aal2} and the
present work is the construction of an operator $Q_\sigma$ from a
subspace of $\mathbf L_2(\mathbb R,dx)$ into  itself (denoted by
$T_m$ in \cite{aal2}; see \eqref{tm}) such that
\begin{equation}
K_\sigma(t,s)=\langle Q_\sigma(1_{[0,t]}),
Q_\sigma(1_{[0,s]})\rangle_{\mathbf L_2(\mathbb
R,dx)}=E[X_\sigma(t)X_\sigma(s)^*].
\label{magic}
\end{equation}
The natural isometry from $\mathbf L_2(\mathbb R,dx)$ into the
white noise space allows then to proceed by doing analysis in the
Gelfand triple associated with the Kondratiev spaces, see Section
\ref{wickito}. Although there are numerous possible other Gelfand
triples, V\r{a}ge's inequality, see \eqref{vage}, is a feature
which seems characteristic of this triple and is very useful in
the computations. In the present section we present some general
results on the existence and properties of the operator
$Q_\sigma$. To develop the associated stochastic integral and
Ito's formula, one has to relate the properties of $\sigma$ and
of $Q$. Not all the arguments go through for general $\sigma$'s.
Details will be presented in forthcoming publications.\\

\begin{Tm}
\label{tmsigma} Let $\sigma$ be a positive measure subject to
\begin{equation}
\int_{\mathbb R}\frac{d\sigma(u)}{1+u^2}<\infty,
\label{dsigma}
\end{equation}
and assume that ${\rm dim}~\mathbf L_2(d\sigma)=\infty$. There
exists a possibly unbounded operator $Q$ from $\mathbf
L_2(d\sigma)$ into $\mathbf L_2(\mathbb R, dx)$, with domain
containing the Schwartz space, such that \eqref{magic} holds:
\begin{equation*}
K_\sigma(t,s)=\langle Q(1_{[0,t]}), Q(1_{[0,s]})\rangle_{\mathbf
L_2(\mathbb R,dx)},
\end{equation*}
and
\[
X_\sigma(t)=\widetilde{Q1_{[0,t]}}.
\]
\end{Tm}

The operator $Q$ in the preceding theorem is not unique.
In the previous cases, a specific choice of $Q$ was made by
recipe, which depended on the special structure of $\mathbf L _2(d\sigma)$.
In general there seems no natural way to chose a specific $Q$.
We have dropped therefore the index $\sigma$ in the notation.\\

{\bf Proof of Theorem \ref{tmsigma}:} We proceed in a number of
step.\\

STEP 1: {\sl Let $W$ be a unitary map from $\mathbf L_2(d\sigma)$
onto ${\mathbf L_2(\mathbb R,dx)}$, and let
\[
f_1=Wb_1,
\]
where $b_1$ denotes the function
\[
b_1(u)=\frac{1}{\sqrt{1+u^2}}\in\mathbf L_2(d\sigma).
\]
Then:
\[
\|f_1\|^2_{\mathbf L_2(\mathbb R,dx)}=\|b_1\|^2_{\mathbf
L_2(d\sigma)}=\int_{\mathbb R}\frac{d\sigma(u)}{u^2+1}<\infty.
\]
}

This is clear from the unitarity of $W$.\\

STEP 2: {\sl The operator $M_u$  of multiplication by the
variable $u$ is {\it a priori} an unbounded operator in $\mathbf
L_2(d\sigma)$. It is densely defined and self-adjoint  in
$\mathbf L_2(d\sigma)$.}\\

This follows from \eqref{eq:new}, which is true for any measure
$\sigma$ satisfying \eqref{dsigma}.\\

It follows from the preceding step that the operator
\[
T=WM_uW^*.
\]
is a self-adjoint operator in $\mathbf L_2(\mathbb R, dx)$. The
operator $Q$
\begin{equation}
\label{QQ}
Q\psi=\sqrt{1+T^2}\widehat{\psi}(T)f_1
\end{equation}
can therefore be computed using the spectral theorem. We claim $Q$
satisfies \eqref{magic} and \eqref{normQ}. This is done in the
next two steps.\\
\

STEP 3: {\it \eqref{normQ} holds.}\\

Indeed, using the spectral theorem we have

\[
\begin{split}
\|Q\psi\|_{\mathbf L(\mathbb R,dx)}^2
&=\|(\sqrt{I+T^2}\widehat{\psi}(T)f_1\|_{\mathbf L(\mathbb R,dx)}^2\\
&=\|W^*\sqrt{I+T^2}\widehat{\psi}(T)Wb_1\|^2_{\mathbf L_2(d\sigma)}\\
&=\int_{\mathbb R
}|\sqrt{1+u^2}\widehat{\psi}(u)\frac{1}{\sqrt{1+u^2}}
|^2d\sigma(u)\\
&=\int_{\mathbb R}|\widehat{\psi}(u)|^2d\sigma(u).
\end{split}
\]

STEP 4: {\it Set}
\[
X_\sigma(t)=\widetilde{Q1_{[0,t]}}.
\]
{\it Then }
\[
E[X_\sigma(t)X_\sigma(s)^*]= \langle Q1_{[0,t]},Q1_{[0,s]}
\rangle_{\mathbf L(\mathbb R,dx)}.
\]
Indeed, the function $1_{[0,t]}$ belongs to the domain of $Q$ and
the claim is a direct consequence of the isometric imbedding of
${\mathbf L(\mathbb R,dx)}$ inside $\mathcal W$.
\mbox{}\qed\mbox{}\\

\begin{Cy}
Let $\sigma$ be a positive measure subject to
\begin{equation}
\label{1plusu}
\int_{\mathbb R}\frac{d\sigma(u)}{1+|u|}<\infty,
\end{equation}
and let $\left\{X_\sigma(t)\right\}$ be the associated process as
in Theorem \ref{tmsigma}. Then, $\left\{X_\sigma(t)\right\}$ has
a continuous version. (By this we  mean \cite{MR2001996} that
$\left\{X_\sigma(t)\right\}$ agrees a.e. with some time-continuous
process.)
\end{Cy}

{\bf Proof:} We will prove this as an application of Kolmogorov's
test for the existence of a continuous version, \cite[p. 14]
{MR2001996}. Because $\left\{X_\sigma(t)\right\}$ is Gaussian,
there exists $K$ independent of $t,s$ such that
\[
E[|X_\sigma(t)-X_\sigma(s)|^4]=K\left(E[|X_\sigma(t)-X_\sigma(s)|^2]\right)^2.
\]
On the other hand,
\[
E[|X_\sigma(t)-X_\sigma(s)|^2]=2{\rm Re}~r(t-s)=2\int_{\mathbb
R}\frac{1-\cos((t-s)u)}{u^2}d\sigma(u).
\]
Using \eqref{1plusu} we now show that
\begin{equation}
{\rm Re}~r(t)\le C|t|,\quad |t|\in[0,1].
\label{kol}
\end{equation}
The result will then follow from Kolmogorov's continuity
criterion (see for instance \cite[Theorem 2.2.3, p. 14]{MR2001996} for the latter).\\

To prove \eqref{kol} we proceed in a way similar as in
\cite{ptrf} as follows. We compute first
\[
\begin{split}
\int_0^1\frac{1-\cos(tu)}{u^2}d\sigma(u)&=2\int_0^1
t^2\frac{\sin\left(\frac{tu}{2}\right)^2}{t^2u^2}d\sigma(u)\\
&\le K_1t^2\quad (K_1>0\quad\mbox{{\rm independent of}}\,\, t)\\
&\le K_1|t|\quad{\rm for}\quad |t|\in[0,1].
\end{split}
\]
Furthermore, using the mean-value theorem for the function
$u\mapsto \cos(tu)$ we have
\[
1-\cos(tu)=t^2u\sin(t\xi_t),\quad \xi_t\in[0,u].
\]
Thus
\[
\begin{split}
\int_1^\infty\frac{1-\cos(tu)}{u^2}d\sigma(u)&=t^2\int_1^\infty\sin(t\xi_t)
\frac{d\sigma(u)}{u}\\
&\le t^2\int_1^\infty\frac{d\sigma(u)}{u}\\
&\le K_2t^2, \quad\mbox{{\rm where we use}}\,\, \eqref{1plusu},\\
&\le K_2|t|\quad{\rm for}\quad |t|\in[0,1],
\end{split}
\]
where $K_2>0$ is independent of $t$. Inequality \eqref{kol}
follows and hence the result.
\mbox{}\qed\mbox{}\\

In Section \ref{secQ} we defined a specific operator $Q_\sigma$
on the Schwartz functions, and then defined $Q1_{[0,t]}$ by
approximation. Here we have used the spectral theorem. Still it
is possible to compute $Q1_{[0,t]}$ via approximating sequences.
We note that, in view of \eqref{dsigma}, the measure
\[
d\mu(u)=\frac{d\sigma(u)}{1+u^2}
\]
satisfies the following property:
\begin{equation}
\label{sigmafinite}
\forall \epsilon>0, \exists \mathbb K \,\, \mbox{{\rm compact and
such that}}\,\, \mu(\mathbb R\setminus \mathbb K)\le \epsilon.
\end{equation}
The arguments in Proposition \ref{10t} can be adapted as follow.
We take as a special sequence
\[
s_n=1_{[0,t]}\star k_{1/n}
\]
where $k_{1/n}$ is defined via \eqref{eq:sketch}. Then
\[
\widehat{s_n}(u)=\chi_t(u)\cdot e^{-\frac{u^2}{n^2}}.
\]
Instead of \eqref{eq_new} we write (for the special sequence
$(s_n)_{n\in\mathbb N}$ at hand)
\[
\begin{split}
\|Q_\sigma s_n-Q_\sigma s_m\|^2_{\mathbf L_2(\mathbb
R,dx)}&=\int_{\mathbb R}
|\widehat{s_n}(u)-\widehat{s_m}(u)|^2d\sigma(u)\\
&=\int_{\mathbb
R}|\chi_t(u)|^2(e^{-\frac{u^2}{n^2}}-e^{-\frac{u^2}{m^2}})^2d\sigma(u)\\
&= \int_{\mathbb K}
|\chi_t(u)|^2(e^{-\frac{u^2}{n^2}}-e^{-\frac{u^2}{m^2}})^2d\sigma(u)+\\
&\hspace{5mm}+
 \int_{\mathbb
R\setminus\mathbb K}
|\chi_t(u)|^2(e^{-\frac{u^2}{n^2}}-e^{-\frac{u^2}{m^2}})^2 d\sigma(u).
\end{split}
\]
where $\mathbb K$ is a compact to be determined.\\

We first focus on the second integral in the last equality above.
Set
\[
\sup_{m,n}|\chi_t(u)|^2(1+u^2)(e^{-\frac{u^2}{n^2}}-e^{-\frac{u^2}{m^2}})^2
=M<\infty,
\]
and recall that $d\mu(u)=\frac{d\sigma(u)}{1+u^2}$. Then
\[
\begin{split}
\int_{\mathbb R\setminus\mathbb K}
|\chi_t(u)|^2(e^{-\frac{u^2}{n^2}}-e^{-\frac{u^2}{m^2}})^2
d\sigma(u)=\\
&\hspace{-2cm}=
 \int_{\mathbb R\setminus\mathbb K}
|\chi_t(u)|^2(1+u^2)(e^{-\frac{u^2}{n^2}}-e^{-\frac{u^2}{m^2}})^2
d\mu(u)\\
&\hspace{-2cm} \le M\mu(\mathbb R\setminus\mathbb K).
\end{split}
\]

For a preassigned $\epsilon>0$, chose now the compact $\mathbb K$
such that
\[
\mu(\mathbb R\setminus\mathbb K)\le \frac{\epsilon^2}{2M}.
\]
Then
\begin{equation}
\label{eq1}
\int_{\mathbb R\setminus\mathbb K}
|\chi_t(u)|^2(e^{-\frac{u^2}{n^2}}-e^{-\frac{u^2}{m^2}})^2
d\sigma(u)\le \frac{\epsilon}{2}.
\end{equation}

Since $K$ is compact there exists $N$ such that, for $n,m$ larger
than $N$, the integral
\begin{equation}
\int_{\mathbb K}
|\chi_t(u)|^2(e^{-\frac{u^2}{n^2}}-e^{-\frac{u^2}{m^2}})^2d\sigma(u)\le
\frac{\epsilon}{2}.
\label{eq2}
\end{equation}
Therefore, for $n,m\ge N$,
\[
\|Q_\sigma s_n-Q_\sigma s_m\|_{\mathbf L_2(\mathbb R,dx)}\le
\epsilon,
\]
and the sequence $(Q_\sigma s_n)_{n\in\mathbb N}$ is a Cauchy
sequence in the norm of $\mathbf L_2(\mathbb R,dx)$. This
provides a constructive way to compute $Q1_{[0,t]}$. To see
that the obtained limit gives the same value as the
one obtained from the spectral theorem, it suffices to let $n$ go to infinity
in \eqref{eq1} and \eqref{eq2}.

\section{Concluding Remarks}
\label{Rem}
\setcounter{equation}{0}
We conclude with comments
comparing our approach with the literature.\\

{\bf 1.} As in \cite{ptrf}, no adaptability of the integrand with
respect to an underlying filtration has been made. In this sense,
one may regard the integral defined here in fact as a
Wick-Skorohod integral.\\

{\bf 2.} Motivated in part by questions in physics, e.g.,
\cite{MR1244577} and \cite{MMP10, OA10}, there has been a recent
increase in the use of operator theory in stochastic processes,
as reflected in e.g. references \cite{ptrf, MR2414165,
aal2, MR1244577,MR2444857,MR1408433}. In addition, we call
attention to the papers \cite{JoSo09a, JoSo09b, HN09, JO09, LS09}
and the papers cited there. In our present approach, we have been
using tools from the cross roads of harmonic analysis and
stochastic process, as are covered in
\cite{DyMc70, MR1408433, ItMc65, Mc69}.

\setcounter{equation}{0}
\bibliographystyle{plain}
\def\cprime{$'$} \def\lfhook#1{\setbox0=\hbox{#1}{\ooalign{\hidewidth
  \lower1.5ex\hbox{'}\hidewidth\crcr\unhbox0}}} \def\cprime{$'$}
  \def\cprime{$'$} \def\cprime{$'$} \def\cprime{$'$} \def\cprime{$'$}

\end{document}